\magnification=1200
\baselineskip=0.20in

%
%
%
%
%

\def\C{{\bf C}}

\def\P{{\bf P}}
\def\Z{{\bf Z}}

\def\OO{{\cal O}}
\def\e{{\epsilon}}

\def\a{\alpha}

\def\d{\delta}
\def\s{\sigma}

\def\l{\lambda}
\def\Z{{\bf Z}}
\def\D{{\Delta}}
\def\d{{\delta}}
\def\U{\cup}
\def\lra{\longrightarrow}

\def\mapright#1{\smash{\mathop{\longrightarrow}\limits^{#1}}}

\def\mapdown#1{\downarrow\rlap{$\vcenter{\hbox{$\scriptstyle#1$}}$}}
\def\mapup#1{\uparrow\rlap{$\vcenter{\hbox{$\scriptstyle#1$}}$}}
\font\mysmall=cmr8 at 8pt    
\centerline{ 
\bf{ON INFINITE DIMENSIONAL GRASSMANNIANS} 
}
\medskip
\centerline{\bf{AND THEIR QUANTUM DEFORMATIONS}}
\bigskip
\centerline{R. FIORESI}
\bigskip
\centerline{Dipartimento di Matematica, Universita' di Bologna}
\centerline{Piazza Porta San Donato 5, 40126 Bologna, Italy}
\centerline{{\mysmall e-mail: fioresi@dm.unibo.it}}
\bigskip
\centerline{C. HACON}
\bigskip
\centerline{Department of Mathematics, UC Riverside}
\centerline{Riverside, CA 92521-0135, USA}
\centerline{{\mysmall e-mail: hacon@math.ucr.edu}}
\bigskip
\centerline{{\bf Abstract}}

{\mysmall An algebraic approach is
developed to define and study infinite dimensional grassmannians.
Using this approach a quantum deformation is obtained for both
the ind-variety union of all finite dimensional
grassmannians $G_{\infty}$, 
and the Sato grassmannian $\widetilde{UGM}$ introduced
by Sato in [Sa1], [Sa2]. They are both quantized as homogeneous
spaces, that is together with a coaction of a quantum infinite
dimensional group. At the end, an infinite dimensional version of
the first theorem of invariant theory is discussed for both the
infinite dimensional special linear group and its quantization.}
\bigskip  
{\bf 1. Introduction}
\medskip
A definition of the infinite dimensional Sato grassmannian 
is first introduced
by Sato in [Sa1], [Sa2], where he explicitly exhibits the points
as infinite dimensional matrices. Sato proves the remarkable fact that the
points of the Sato grassmannian 
$\widetilde{UGM}$ are in one to one correspondence with
the solutions of the KP hierarchy. 

A few years later Segal and Wilson [SW], using mainly analytic techniques,
explore more deeply this correspondence. 

In a later work [PS] 
Pressley and Segal study more extensively, along the same lines,
an infinite dimensional grassmannian closely related to
the Sato Grassmannian. In particular they give a stratification and 
a Plucker embedding of it. They also produce an action of a certain
infinite dimensional linear group realizing it as an
infinite dimensional homogeneous space.
Though their definition appears quite different from Sato's one, they
essentially describe the same geometrical object, but in a slightly
more general setting.

A more geometrical approach to the same subject is taken by
Mulase in [Mu1], [Mu2]. He constructs the Sato grassmannian
as a scheme of which
he gives the functor of points. Also Plaza-Martin [PM] take the
same approach, with special attention given to the physical applications.

Together with the Sato grassmannian, Sato, as well as all the above
mentioned authors, introduces what we denote by $G_{\infty}$, the
union of all finite dimensional grassmannians. $G_{\infty}$ turns out to be
an ind-variety [Ku] and it is dense in various topologies inside
$\widetilde{UGM}$. $G_{\infty}$ is an interesting object in itself. 
Using the points of $G_{\infty}$ expressed as infinite wedge products,
in [Ka] Kac constructs an infinite dimensional representation
of an infinite dimensional general linear group and shows
the correspondence between points of the infinite dimensional
grassmannian and solutions
of the KP hierarchy with algebraic methods.

In the present work we want to study the infinite dimensional grassmannians 
$\widetilde{UGM}$ and  $G_{\infty}$ using
only algebraic methods. This approach turns to be the most
natural for our goal,
that is to obtain their quantum deformations.   
We will include proofs of statements about grassmannians
generally known in the literature
whenever an appropriate reference is not available.

This paper is divided in three parts. 

In the first part, \S 2, we consider the inverse and direct limit of  
the coordinate rings $k[\d_{m,n}]$ of the finite dimensional grassmannian over
the algebraically closed field $k$. Then we give an explicit presentation
for the inverse limit $\widehat{k[\d_{\infty}]}$ and the direct limit
$k[\d_{\infty}]$. 
We also prove that $\widehat{k[\d_{\infty}]}$ and ${k[\d_{\infty}]}$
can be in some sense regarded as the homogeneous coordinate rings of 
$G_{\infty}$ and $\widetilde{UGM}$. In fact the closed points
of Proj($\widehat{k[\d_{\infty}]}$) and of Proj($k[\d_{\infty}]$)
turn to be in one-to-one correspondence with the
points of $G_{\infty}$ and $\widetilde{UGM}$ respectively.
Both  $G_{\infty}$ and
$\widetilde{UGM}$ admit an action of the infinite dimensional
special linear group $SL_{\infty}$ given by the union of all finite
dimensional special linear groups over $k$.
We also show that there is a corresponding coaction of the homogeneous
coordinate
ring of the ind-variety $SL_{\infty}$ on both 
$\widehat{k[\d_{\infty}]}$ and ${k[\d_{\infty}]}$.

In the second part of the paper, \S 3, we repeat these same constructions
in the quantum groups setting. We give explicit quantum deformations
for both the ind-variety $G_{\infty}$ and the Sato grassmannian
$\widetilde{UGM}$.  
Proceeding in the same way as in \S 2, we take the inverse and
direct limit of the quantum finite dimensional grassmannian
${k_q[\D_{m,n}]}$ ([Fi1], [TT]). We obtain two non commutative
rings, $\widehat{k_q[\D_{\infty}]}$ and ${k_q[\D_{\infty}]}$
deformations of 
$\widehat{k[\d_{\infty}]}$ and ${k[\d_{\infty}]}$ respectively
that we call quantum $G_{\infty}$ and 
quantum Sato grassmannian $\widetilde{UGM}$.
We give an explicit presentation for both of them.
$G_{\infty}$ and
$\widetilde{UGM}$ are quantized as homogeneous spaces, that is
there is a well defined coaction of the quantum special linear 
infinite dimensional group $k_q[SL_{\infty}]$ on them.

In the last part, \S 4, we examine the following problem of classical
invariant theory for the infinite dimensional case: given the natural
right action of the special linear group of order $r$, 
$SL_{r,0}(k)$ on the matrix algebra, find the
$SL_{r,0}(k)$-invariants.
In complete analogy to what happens in the finite dimensional case,
the ring of invariants in the infinite dimensional case coincides
with $\widehat{k[\d_{\infty}]}$ the homogeneous coordinate ring for the
ind-variety $G_{\infty}$.
We then obtain the corresponding results for the quantum case,
generalizing the results in the paper [FH]. 

The first author wishes to thank Prof. V.S. Varadarajan 
and Prof. I. Dimitrov for many fruitful
discussions and Prof. R. Achilles for helpful comments.
\bigskip
{\bf 2. The infinite dimensional grassmannians $G_{\infty}$ and 
$\widetilde{UGM}$}
\medskip
Let $k$ be an algebraically closed field of characteristic 0.

Let $G_{(m,n)}$ be the grassmannian of $m$ subspaces in a vector space of
dimension $N=m+n$.   

An element of $G_{(m,n)}$ is represented by a
$N \times m$ matrix. We will assume (following Sato [Sa2]) that the
row indices go from $-m$ to $n-1$  while
the column indices go from $-m$ to $-1$.

Let $k[a_{i,j}]_{m,n}$ be the coordinate ring of the algebra of
the $N \times N$ matrices, where we assume that both row and column
indices go from $-m$ to $n-1$.

The homogeneous coordinate ring of $G_{(m,n)}$ is isomorphic to the subring of
the matrix ring $k[a_{i,j}]_{m,n}$ generated by
the determinants $d_{l_{0} \dots l_{m-1}}$ of the minors
obtained by taking the columns $-m \dots -1$ and the rows
$l_{0} \dots l_{m-1}$. We will denote such subring by
$k[\d_{m,n}]$ and the above mentioned set of determinants by
$\d_{m,n}$.


{\bf Definition (2.1)}. Let $m'\geq m$, $n'\geq n$.
Define the inverse family of rings:
$$
\matrix
{ 
k[\d_{m',n'}] & \mapright {e_{(m',n',m,n)}} & k[\d_{m,n}] \cr
}
$$
$$
e_{(m',n',m,n)}(d_{j_0 \dots j_{m'-1}})=\cases{
d_{l_0 \dots l_{m-1}} & if $(j_0 \dots j_{m'-1})=
(-m' \dots -m-1,l_0 \dots l_{m-1})$ 
\cr
& and $-m \leq l_0 < ... < l_{m-1} \leq n-1$ 
\cr \cr
0 & otherwise
} 
$$

with $-m'\leq j_0<\dots <j_{m'-1} \leq n'-1$.

We define:
$$
\widehat {{k}[\d_{\infty}]}=\lim_{\longleftarrow} k[\d_{m,n}],
$$
and denote the induced maps by
$$
e_{(m,n)}: \widehat {k[\d_{\infty}]} \lra k[\d_{m,n}].
$$
We observe that the maps $e_{(m',n',m,n)}$ are induced by maps
$$E_{(m',n',m,n)}:k[a_{i,j}]_{m',n'}\lra k[a_{i,j}]_{m,n}$$  
defined by
$$E_{(m',n',m,n)}(a_{i,j})=a_{i,j},\ \forall \ -m \leq i\leq n-1,\
-m\leq j\leq -1$$
$$E_{(m',n',m,n)}(a_{i,j})=1,\ \forall \ -m'\leq i=j\leq -m-1$$
$$E_{(m',n',m,n)}(a_{i,j})=0\ \hbox{otherwise.}$$
  
We define:
$$
k[M_{\infty }]=\lim_{\longleftarrow} k[a_{i,j}]_{m,n}.
$$

{\bf Remark (2.2)}. Any element $b\in k[M_{\infty }]$
is an element of the form $b=\{ b_{(m,n)}\}$ such that
$b_{(m,n)}\in k[a_{i,j}]_{m,n}$ and for all $m'\geq m$, $n'\geq n$,
$E_{(m',n',m,n)}(b_{(m',n')})=b_{(m,n)}$.
Similarly any element  $x \in \widehat {k[\d_{\infty}]}$ is of the form
$x=\{ x_{(m,n)}\}$ such that
$x_{(m,n)}\in k[\d _{i,j}]_{m,n}$ and for all $m'\geq m$, $n'\geq n$,
$e_{(m',n',m,n)}(x_{(m',n')})=x_{(m,n)}$.
\medskip
There is a corresponding direct system of
inclusions of projective
varieties
$$G_{(m,n)}\lra G_{(m',n')}\ \ \forall\ \ m'\geq m,\ n'\geq n.$$
Define $$G_{\infty }=\lim_{\longrightarrow}G(m,n).$$ Notice that
$G_{\infty }=\cup G_{(m,n)}$.
\medskip
We want to view $G_{\infty }$ as a projective
ind-variety. 
\medskip
{\bf Definition (2.3)}.
An {\it ind-variety over $k$} is a set $X$ together with a filtration:
$$
X_0 \subset X_1 \subset X_2 \subset ...
$$
such that

1) $\U_{n \geq 0} X_n = X$

2) Each $X_n$ is a finite dimensional variety over $k$ such that
the inclusion $X_n \subset X_{n+1}$ is a closed immersion. 

(See [Ku] for more details).
\medskip
The ind-variety $X$ is naturally a topological space, $U\subset X$ being {\it
open} in $X$ if and only if, for each $n$, $U\cap X_n$ is open in $X_n$.
The sheaf of {\it regular functions on $X$} is defined by $\OO _X := 
\matrix{\lim_{} \cr
{\leftarrow}} \OO _{X_n}$.
$X$ is said to be a {\it locally
projective ind-variety} if it admits a filtration such that
each $X_n$ is projective.
We will say that $X$ is a {\it
projective ind-variety} if it admits a filtration 
$X_n$ and a line bundle $L$ such that each restriction
$L|X_n$ is very ample and the corresponding maps
$$H^0(X_n,L|X_n)\longrightarrow H^0(X_{n-1},L|X_{n-1})$$
are surjective. In other words, for each $n$ there are compatible
closed immersions $X_n \hookrightarrow \P ^{N_n}=\P (H^0(X_n,L|X_n)^\vee)$
with coordinate rings generated by
$H^0(X_n,L|X_n)$ and hence a closed immersion of
ind-varieties $X\hookrightarrow \P ^\infty=\cup \P (H^0(X_n,L|X_n)^\vee)$.
We define 
$$H^0(X,L):=\lim_{\longleftarrow} H^0(X_n,L|X_n).$$
Let $S(\P^N)=\oplus _{d\geq 0}H^0(\P ^N, \OO _{\P ^N}(d))$ be
the homogeneous coordinate ring of $\P ^N$ and $I(X_n)$ be the homogeneous
ideal of $X_n\subset \P ^{N_n}$, then the homogeneous coordinate ring
of $X_n\subset \P ^{N_n}$ is given by $S(X_n)=S(\P ^N)/I(X_n)$.
We define the {\it homogeneous coordinate ring} of the projective ind-variety
$X\subset \P ^\infty$ to be 
$$S(X):=\lim_{\longleftarrow}S(X_n).$$
{\bf Theorem (2.4)}. {\it $G_{\infty }$ is a projective ind-variety, with
homogeneous coordinate ring $\widehat {k[\d_{\infty}]}$.}

{\bf Proof}.
It is well known that the maps $G_{(m,n)}\lra G_{(m',n')}$ are closed
immersions (when defined). Let $X_n:=G_
{(n,n)}$, then we have closed immersions $G_{(m,n)}\hookrightarrow
X_{(n+m,n+m)}$. Therefore $X=\cup X_n = \cup G_{(m,n)}=G_{\infty }$
is an ind-variety.
For each $n>0$ we have 
the Pl\"ucker embeddings
$X_n=G_{n,n}\hookrightarrow \P (\wedge ^n \C ^{2n})=\P ^{N_n}$.
The homogeneous coordinate ring $S(\P ^{N_n})$
is generated by elements $x_I$ where $I=\{i_1,...,i_n\}$ such that
$-n\leq i_1 <i_2 < ... < i_n \leq n-1$ and for $n'\geq n$, 
the closed immersions
$\P ^{N_n}\hookrightarrow \P ^{N_{n'}}$ correspond to the 
surjective homomorphisms
$$S(\P ^{N_{n'}})
\lra S(\P ^{N_n})$$
defined by
$$
\e _{(m',n',m,n)}(x_{i_1 \dots i_{n'}})=\cases{
x_{i_1 \dots i_{n}} & if $(i_1 \dots i_{n'})=
(-n' \dots -n-1,i_1 \dots i_{n})$ 
\cr
& and $-n \leq i_1 < ... < i_{n} \leq n-1$ 
\cr \cr
0 & otherwise
} 
$$
The homogeneous coordinate ring of the projective
ind-variety $\P^\infty=\cup _{n>0}\P^{N_n}$ is just generated by
$\matrix{\lim_{}\cr\longleftarrow} H^0(\P ^{N_n}, \OO _{\P ^{N_n}}(1))$.
The line bundle $L|X_n$ is just the pull back of $\OO _{
\P (\wedge ^n \C ^{2n})}(1)$.
The immersions $\P^{N_n}\lra \P ^{N_{n'}}$ and $X_n\lra X_{n'}$ are compatible .
The corresponding
homogeneous coordinate ring is $S(X_n)=k[\d_{n,2n}]$.
The maps induced by the inclusions
$X_n\lra X_{n'}$ are just the maps $e_{(n',n',n,n)}:k[\d_{n',2n'}]
\longrightarrow k[\d_{n,2n}]$. 
Therefore the homogeneous coordinate ring of $X$ is given by
$$\lim_{\longleftarrow}k[\d_{n,2n}]=\lim_{\longleftarrow}k[\d_{m,n}]=
\widehat {k[\d_{\infty}]}.$$
QED.
\medskip

We now turn our attention to the Sato grassmannian $\widetilde{UGM}$
and its relation with $G_{\infty}$.

{\bf Definition (2.5)}. Let $m'\geq m$, $n'\geq n$.
Define a direct family of rings:
$$
\matrix
{
k[\d_{m,n}] & \mapright{r_{(m,n,m',n')}} & k[\d_{m',n'}] \cr
d_{l_0 \dots l_{m-1}} & \longmapsto &
d_{-m', \dots -m-1,l_0 \dots l_{m-1}}
}
$$
for $-m \leq l_0 < \dots < l_{m-1} \leq n-1$.

It is easy to see using the Pl\"ucker relations that this map is well defined.
Moreover,
the map $e_{(m,n,m',n')}$ is a left inverse for $r_{(m,n,m',n')}$ and 
in particular $r_{(m,n,m',n')}$ is injective.
We define: 
$$
k[\d_{\infty}]=\lim_{\longrightarrow} k[\d_{m,n}].
$$
Denote with $r_{(m,n)}$ the induced inclusions 
$k[\d_{m,n}]\lra k[\d_{\infty}]$.

{\bf Definition (2.6)}. Define Maya diagram of virtual cardinality 0
(or shortly a {\it Maya diagram}) 
a strictly increasing sequence $a_{\bullet}=
\{ a_i \}$, $i\geq 1$, such that $a_i\in \Z$
and $a_i=i$ for all $i>>0$. Define the order $|| a_{\bullet}||$ of
a Maya diagram to be the smallest number $i$ such that $a_j=j$
for all $j\geq i$. Any sequence $l_*=l_1,\dots ,l_{m}$ with
$l_1 \leq ... \leq l_{m}\leq m$ induces
a Maya diagram $a_{\bullet}=\tilde {l_*}$ of order at most $m+1$ defined by
$a_i=l_{i}$ for all $1\leq i\leq m$ and $a_i=i$ for all $i\geq m+1$. 
For any Maya diagram $a_{\bullet}$, let $a_{\leq m}$
denote the ordered set $a_1<\dots <a_m$. Clearly
if $||a_{\bullet}||\leq m$, then $a_{\bullet}=\widetilde {a_{\leq m}}$.

Given a Maya diagram $a_{\bullet}$ of order $m+1$
with $|a_1| \geq -n+1$, we wish to define
corresponding elements $d_{a_{\bullet}}\in {k}[\d_{\infty}]$, and
$\widehat {d} _{a_{\bullet}}\in \widehat {k[\d_{\infty}]}$.
Define
$$
d_{a_{\bullet}}=r_{(m,n)}d_{a_{\leq m}}
$$
where $a_{\bullet}=\widetilde {a_{\leq m}}$.

$k[\delta_\infty ]$ is generated as a ring by the $d_{a_\bullet}$, since
it is generated by the images of $k[\delta_{m,n}]$ under
$r_{(m,n)}$.
 
We define a map
$$
\matrix{
\matrix{
\rho _{(m,n)}:k[\d_{\infty}] & \longrightarrow & k[\d_{m,n}]
} \cr \cr
\rho _{(m,n)} (d_{a_{\bullet}})= 
\cases{
d_{a_{\leq m}} &  for all $m\geq ||a_{\bullet}||$, $n\geq |a_1|$ \cr \cr
0 &  otherwise
}
\cr
}
$$ 
Then we define
$$
\widehat {d}_{a_{\bullet}}=\{ \rho _{(m,n)}
d_{a_{\bullet}} \} \in \widehat{k[\d_{\infty}]}.
$$

{\bf Proposition (2.7)}. 

{\it a) There is an injection 
$I:{k}[\d_{\infty}]\lra \widehat {k[\d_{\infty}]}$
which sends $d _{a_{\bullet}}$ to $\widehat {d}_{a_{\bullet}}$.

b) The image of $I$ is dense in $\widehat {k[\delta _\infty ]}$
for the inverse limit topology on $\widehat {k[\delta _\infty ]}$
induced by the discrete topology on each $k[\delta _{m,n}]$.
}

{\bf Proof}. (a) It suffices to check that for $m'\geq m$ and $n'\geq n$,
one has $$e_{(m',m,n',n)}(\rho _{(m',n')}d_{a_\bullet })=
\rho _({m,n)}d_{a_\bullet }.$$

\noindent
(b) Define a fundamental set of neighbourhoods of $0 \in\widehat {
 k[\d_{\infty}]}$
by $U_0=\widehat {k[\d_{\infty}]}$, $U_{k}=e_{(k,k)}^{-1}(0)$. 
For any $x=\{ x_{(m,n)} \}
\in \widehat {k[\d_{\infty}]}$, we must define a Cauchy sequence 
$\{y_k\}=\{\sum_{i=1}^{m_k} b^k_i\widehat{d}_{a_{\bullet}^{i,1}}
\dots \widehat{d}_{a_{\bullet}^{i,j_i}}\}$ lying in the image of $I$ and
converging to $x$.


Since $e_{(k,k)}(x)
\in k[\d_{k,k}]$ we have 
$$e_{(k,k)}(x)=\sum_{i=1}^{m_k} b^k_i {d}_{L_k^{i,1}} \dots {d}_{L_k^{i,j_i}}.
$$
where $L_k^{i,j}=(l_{k,1}^{i,j}...l_{k,m}^{i,j})$,  
$-k \leq l_{k,1}^{i,j}<\dots <l_{k,m} ^{i,j} \leq k-1$. 
Set 
$$
y_k=\sum_{i=1}^{m_k} b^k_i {\widehat{d}}_{\widetilde{L_k^{i,1}}} \dots \widehat
{d}_{\widetilde{L_k^{i,j_i}}}.
$$
Notice that the summation is finite and hence $y_k$ is in the image of $I$.
For any $k'>k$ we have, $e_{(k,k)}(y_{k'})=e_{(k,k)}(x)$
i.e. $y_{k'}-x \in U_k$.
Hence $\{y_k\}$ is a Cauchy sequence converging to $x$.
QED
\medskip

We now give a presentation of the rings ${k}[\d_{\infty}]$ and
$\widehat {k[\d_{\infty}]}$. 

Define $k[\xi _{a_{\bullet }}]$ to be the ring generated by the independent
variables $\xi _{a_{\bullet }}$, where ${a_{\bullet }}$ is any
Maya diagram of virtual cardinality 0.
 
There is a natural map $\phi :k[\xi _{a_{\bullet }}]  \lra  k[\d _{\infty }]$
such that $\xi _{a_{\bullet }}  \lra  d _{a_{\bullet }}$.
This induces a topology on $k[\xi _{a_{\bullet }}]$ for which a fundamental set
of neighborhoods is given by $V_k:=\phi ^{-1} I^{-1}U_k$. Let
$\widehat {k[\xi _{a_{\bullet }}]}$ be the completion of $
k[\xi _{a_{\bullet }}]$ with respect to the above
topology. In particular the elements of $\widehat {k[\xi _{a_{\bullet }}]}$
are of the form
$\sum b_i \xi _{a_{\bullet }^{i,1}}\cdots
\xi _{a_{\bullet }^{i,k_i}}$ where $b_i\in k$ and  the
${a_{\bullet}^{i,k}}$ are
any Maya diagrams of virtual cardinality 0.
The corresponding natural map between completions 
$\widehat {\phi } :\widehat {k[\xi _{a_{\bullet }}]} \lra 
\widehat {k[\d _{\infty }]}$ is defined by 
$\xi _{a_{\bullet }}   \lra \widehat{d}_{a_{\bullet }}$.

Define
$P_{(m,n)}$ to be the ideal of Plucker relations in $k[\d_{m,n}]$.
Let $P=\cup P_{(m,n)}$ be the corresponding ideal in $k[\d_{\infty}]$,
and similarly $\widehat {P}=
\matrix{\lim \cr {\longleftarrow}} P_{(m,n)}$
be the corresponding ideal in $\widehat {k[\d_{\infty}]}$.
\vskip1ex
{\bf Theorem (2.8)}. {\it We have ring isomorphisms
$$
\matrix{
i) \quad 
\widehat {k[\d _{\infty }]}\cong \widehat {k[\xi _{a_{\bullet }}]}
/\widehat {P} \cr \cr
ii) \quad k[\d _{\infty }]\cong k[\xi _{a_{\bullet }}]/P 
}
$$
}

{\bf Proof}. (i) The direct limit is an exact functor.

(ii) The inverse limit functor is left exact, and since the inverse system  
$P_{(m,n)}$ is a surjective system, the corresponding inverse
system sequence
$$
0 \lra \widehat {P}\lra \widehat {k[\xi _{a_{\bullet }}]}\lra
\widehat {k[\d _{\infty }]}\lra 0
$$
is also exact.  QED.
\medskip
We want now to relate our constructions with the Sato grassmannian.
In [Sa2]  Sato defines a set of points in an infinite dimensional
projective space. 
Already theorem (2.8) suggests to view $G_{\infty}$
as the set of zeros of the ideal $\widehat{P}$ in an infinite dimensional 
projective space whose coordinate ring is given by
$\widehat{k[\xi _{a_{\bullet }}]}$. 
We want to make this euristic notion more precise and 
to relate the ring ${k[\d_{\infty}]}$ with the Sato
grassmannian.

Assume that the field $k$ has cardinality strictly greater than $\aleph_0$.

Consider the directed system given by rings $R_n=k[z_1,...,z_n]$ 
and homomorphisms of
$k[z_1,...,z_{n'}] \lra k[z_1,...,z_n]$ (for $n'>n$) defined by
$z_i \lra z_i$ for all $i\leq n$, and $z_i\lra 0$ for $i>n$.
This corresponds to an affine ind-variety $A^\infty=\cup  A^n$
given by the inclusion of affine planes $A^n =Spec (R_n) \hookrightarrow
A^{n+1} =Spec (R_{n+1})$. 
Let 
$$
\hat {R}=\lim _{\leftarrow } R_n.
$$
{\bf Lemma (2.9)}.
{\it The set of closed points of $A^\infty$ is in one to one
correspondence with $Spec _m(\hat {R})$.}

{\bf Proof}.
Notice that each $R_n$ injects in a natural way in $ \hat {R}$
and set $R=\cup R_n \subset \hat{R}$.
If $m \subset R$ is maximal, $R/m=E$ is a field. By [La] we have $E=k$
and therefore $m$ is generated by $z_i-k_i$ where $k_i$ 
is the image of $z_i \in R/m$. If $m'$ is a maximal ideal of $\hat {R}$,
and $f:\hat {R}\lra \hat {R}/m'$, then by the previous observations, the 
induced maps $f_i:R_i\lra \hat {R}/m'$ have image contained in $k$.
By the universal property of inverse limits then also $f:\hat {R}\lra k$
is determined by $k_i=f(z_i)$. It is clear that in order for $f$ to be defined,
one must have $k_i=0$ for all but finitely many $i$. QED.

Recall that $\widehat{k[\xi _{a_{\bullet }}]}=
\matrix{\lim_{}\cr\longleftarrow} S(\P ^{N_n})$.
It follows that $\P ^\infty = Proj (\widehat{k[\xi _{a_{\bullet }}]})$
and that the closed points of $\P ^\infty$ are given by 
sequences $k_{a_{\bullet}}\in k$ where ${a_{\bullet}}=0$ for all
but finitely many Maya diagrams of virtual cardinality 0 we have 
${a_{\bullet}}=0$ and two sequences are considered equivalent
if there exists $\lambda \in k^*$ such that $k_{a_{\bullet}}=\lambda k^\prime
_{a_{\bullet}}$.
The assertion can be verified locally on the open cover
$Spec (\widehat{k[\xi _{a_{\bullet }}]}_{(\xi _{\bar {a}_{\bullet }})})$
where $\widehat{k[\xi _{a_{\bullet }}]}_{(\xi _{\bar {a}_{\bullet }})}$
denotes the subring of elements of degree 0 in the localized ring
$\widehat{k[\xi _{a_{\bullet }}]}_{\xi _{\bar {a}_{\bullet }}}$.
The computation is now analogous to the one above for $\hat{R}$.
Similarly one has that the closed points of 
$Proj ({k[\xi _{a_{\bullet }}]})$ correspond to all sequences (not necessarily
bounded) $k_{a_{\bullet}}\in k$ where ${a_{\bullet}}$ runs over all
Maya diagrams of virtual cardinality 0 
and two sequences are considered equivalent
if there exists $\lambda \in k^*$ such that $k_{a_{\bullet}}=\lambda k^\prime
_{a_{\bullet}}$.
It follows that

{\bf Proposition (2.10)}.
{\it
Assume that the cardinality of $k$ is strictly greater than $\aleph_0$.
 
i) The set of closed points of $G_{\infty}$ is in one to one correspondence
with the set of closed points of $Proj (\widehat{k[\d_{\infty}]})$, i.e. with
the sequences 
$\{k_{a_{\bullet}} \in k\}$ satisfying all Pl\"ucker relations,
where ${a_{\bullet}}$  Maya diagram of
virtual cardinality 0 and $k_{a_{\bullet}} =0$ for all but finitely many
Maya diagrams.

ii) The set of closed points of $Proj ({k[\d_{\infty}]})$
is in one to one correspondence with 
the sequences 
$\{k_{a_{\bullet}} \in k\}$ satisfying all Pl\"ucker relations,
where ${a_{\bullet}}$  Maya diagram of
virtual cardinality 0. Moreover we have that those points
coincide with $\widetilde {UGM}$ the Sato grassmannian
(as defined by Sato, [Sa1], [Sa2]).

}
\medskip

{\bf Remark (2.11)}. 
Proposition (2.10) shows that the ring ${k[\d_{\infty}]}$
can be regarded as the ``coordinate ring'' for the Sato grassmannian in
the sense that its maximal ideals are in one-to-one correspondence
with the points of  $\widetilde{UGM}$. Theorem (2.9) allows us
to interpret the Sato grassmannian as the set of closed points 
in an infinite dimensional projective space that are subjected to
the relations $P$. 

\medskip
We now want to define an infinite dimensional special linear group
and show that it has an action on both $G_{\infty}$ and
$\widetilde{UGM}$.

{\bf Definition (2.12)}.
For all $m,n$ positive integers, define $SL_{(m,n)}(k) \cong SL_N(k)$ 
as $N\times N$ matrices with determinant 1, whose row and column indicies are 
between $-m$ and $n-1$.
The inclusions $\psi_{(m,n,m',n')}:SL_{m,n}(k) \lra SL_{m',n'}(k)$ are defined 
for all
$m'\geq m$, $n'\geq n$, $\psi_{(m,n,m',n')}(g)=diag(Id_{m-m'},g,Id_{n-n'})$. 
It is clear that we have an action of $SL_{m,n}(k)$ on $G_{(m,n)}$ for
all $m,n$. We have a corresponding projective system of coordinate
rings:
$$
\matrix{
k[SL_{m',n'}] \quad \mapright{\phi_{(m',n',m,n)}}  \quad k[SL_{m,n}]
\cr \cr
\phi_{(m',n',m,n)}(g_{ij})=
\cases{
g_{ij} & if $-m \leq i,j \leq n-1$
\cr \cr
1  & if $-m' \leq i=j \leq -m-1$, $n \leq i=j \leq n'-1$ 
\cr \cr
0  & otherwise
}
}
$$

{\bf Observation (2.13)}.
$$
SL_{\infty}(k)=_{def} 
\lim_{\longrightarrow} SL_{m,n}=
\U SL_{m,n}(k)$$ 
is an ind-variety with coordinate ring
$$
k[SL_{\infty}]=_{def}\lim_{\longleftarrow} k[SL_{m,n}].
$$
We want now to show that $k[SL_{\infty}]$ has an Hopf algebra structure.
Notice that while in the finite dimensional case this is an immediate
consequence of the fact that the variety $SL_{m,n}(k)$ is a group,
in the infinite dimensional case we need to check the commutativity
of certain diagrams.

{\bf Proposition (2.14)}. {\it 
Let $(m',n')>(m,n)$. The following are commutative diagrams:

i)
$$
\matrix{
k[SL_{m',n'}] & \mapright{\phi_{(m',n',m,n)}} & k[SL_{m,n}] \cr \cr
\mapdown{\Delta_{(m',n')}} & & \mapdown{\Delta_{(m,n)}} \cr \cr
k[SL_{m',n'}] \otimes k[SL_{m',n'}] &
\mapright{\phi_{(m',n',m,n)} \otimes \phi_{(m',n',m,n)}} &
k[SL_{m,n}] \otimes k[SL_{m,n}]
}
$$
where $\Delta_{(m,n)}$ is the comultiplication
in the Hopf algebra $k[SL_{m,n}]$.

ii)
$$
\matrix{
k[SL_{m',n'}] & \mapright{\phi_{(m',n',m,n)}} & k[SL_{m,n}] \cr\cr
\mapdown{\epsilon_{(m',n')}} & & \mapdown{\epsilon_{(m,n)}} \cr\cr
k & \mapright{id} & k
}
$$
where $\epsilon_{(m,n)}$ is the counit in the Hopf algebra $k[SL_{m,n}]$.

iii) 
$$
\matrix{
k[SL_{m',n'}] & \mapright{\phi_{(m',n',m,n)}} & k[SL_{m,n}] \cr \cr
\mapdown{S_{(m',n')}} & & \mapdown{S_{(m,n)}} \cr \cr
k[SL_{m',n'}] & \mapright{\phi_{(m',n',m,n)}} & k[SL_{m,n}]
}
$$
where $S_{(m,n)}$ is the antipode in the Hopf algebra $k[SL_{m,n}]$.
}

{\bf Proof}. Direct check.

{\bf Corollary (2.15)}. {\it $k[SL_{\infty}]$ has an Hopf algebra
structure given by:

a) comultiplication
$$
\matrix{
k[SL_{\infty}] & \mapright{\Delta_{\infty}} & 
k[SL_{\infty}] \hat \otimes k[SL_{\infty}] \cr \cr
\{a_{(m,n)}\} & \mapsto & \{\Delta_{(m,n)}(a_{(m,n)})\}
}
$$

b) counit
$$
\matrix{
k[SL_{\infty}] &  \mapright{\epsilon_{\infty}} & k \cr \cr
\{a_{(m,n)}\} & \mapsto & \epsilon_{(m,n)}(a_{(m,n)})
}
$$

c) antipode
$$
\matrix{
k[SL_{\infty}] & \mapright{S_{\infty}} & k[SL_{\infty}]  \cr \cr
\{a_{(m,n)}\} & \mapsto & \{S_{(m,n)}(a_{(m,n)})\}
}
$$
where $\hat \otimes$ denotes the completed tensor product 
and is given by 
$$
k[SL_{\infty}] \hat \otimes k[SL_{\infty}]=\lim_{\longleftarrow}
k[SL_{m,n}]  \otimes k[SL_{m,n}]
$$
(see [Ku] for more details).}

 
{\bf Proof.} (a) is immediate from proposition (2.14) and from [Ku].
(b), (c) are immediate from proposition (2.14).

\medskip
The group 
$SL_{\infty}(k)$ has an action on both $G_{\infty}$ and $\widetilde{UGM}$.
In order to obtain a quantization of these actions we need to
describe the corresponding coactions of ${k[SL_{\infty}]}$ on 
$\widehat{k[\d_{\infty}]}$ and ${k[\d_{\infty}]}$.

{\bf Observation (2.16)}.
Since $SL_{m,n}(k)$ acts on $G_{(m,n)}$ we have the coaction:
$$
\matrix{
k[\d_{m,n}] & \mapright{\lambda_{(m,n)}} & k[SL_{m,n}] \otimes k[\d_{m,n}] 
\cr \cr
d_{l_0 \dots l_{m-1}} & \longmapsto &
\sum_{m \leq k_0 ... k_{m-1} \leq n-1} g_{l_0k_0}
...g_{l_{m-1}k_{m-1}} \otimes d_{k_0 ... k_{m-1}}
}
$$
One can check the commutativity of the following diagram, for
$m'\geq m$, $n'\geq n$:
$$
\matrix{
k[\d_{m,n}] & \mapright{\lambda_{(m,n)}} & k[SL_{m,n}] \otimes k[\d_{m,n}]
\cr \cr
\mapup{ e_{(m',n',m,n)} } & & \mapup{ \phi_{(m',n',m,n)} \otimes e_{(m',n',m,n)} }
\cr \cr
k[\d_{m',n'}] & \mapright{ \lambda_{(m',n')} } & k[SL_{m',n'}] 
\otimes k[\d_{m',n'}]
}
$$

{\bf Proposition (2.17)}. {\it There is an coaction of ${k[SL_{\infty}]}$ on
$\widehat{k[\d_{\infty}]}$ and on $k[\d_{\infty}]$.}

{\bf Proof}. 
Fix $(m_0,n_0)$. Let $m'> m > m_0$, $n'> n>n_0$.
We have a commutative diagram (see observation (2.16)):
$$
\matrix{
k[\d_{m,n}] & \mapright{\lambda_{(m,n,m_0,n_0)}} & k[SL_{m_0,n_0}] \otimes k[\d_{m,n}]
\cr \cr
\mapup{ e_{(m',n',m,n)} } & & \mapup{ id \otimes
e_{(m',n',m,n)} }
\cr \cr
k[\d_{m',n'}] & \mapright{ \lambda_{(m',n',m_0,n_0)} } & k[SL_{m_0,n_0}]
\otimes k[\d_{m',n'}]
}
$$
$$
\lambda_{(m,n,m_0,n_0)}(d_{l_0...l_{m-1}})=_{def}
\sum_{m \leq k_0 ... k_{m-1} \leq n-1} \tilde g_{l_0k_0}
...\tilde g_{l_{m-1}k_{m-1}} \otimes d_{k_0 ... k_{m-1}}
$$
where:
$$
\tilde g_{ij}=\cases{
g_{ij} & if $-m \leq i,j \leq n-1$
\cr \cr
1  & if $-m' \leq i=j \leq -m-1$, $n \leq i=j \leq n'-1$
\cr \cr
0  & otherwise
}
$$
Hence there is a map:
$$
\matrix{
\widehat{k[\d_{\infty}]} & \mapright{} & 
k[SL_{m_0,n_0}] \otimes \widehat{k[\d_{\infty}]}
}
$$
Now taking the inverse limit
of $k[SL_{m_0,n_0}]$ on the right side we obtain a coaction
of $k[SL_{\infty}]$ on $\widehat{k[\d_{\infty}]}$, that is a
map:
$$
\matrix{
\widehat{k[\d_{\infty}]} & \mapright{} &
k[SL_{\infty}] \hat \otimes \widehat{k[\d_{\infty}]}
}
$$
where $\hat \otimes$ denotes the completed tensor product (see [Ku]).
 
By theorem (2.7) (iii) $k[\d_{\infty}]$ can be identified 
with a subring of $\widehat{k[\d_{\infty}]}$ and one can check
that the given coaction is a well defined coaction 
when restricted to
this subring of $\widehat{k[\d_{\infty}]}$.
QED.
\bigskip
{\bf 3. The quantum infinite dimensional grassmannians $k_q[G_{\infty}]$ and
$k_q[\widetilde{UGM}]$}
\medskip
We want to obtain a deformation of $G_{\infty}$ and $\widetilde{UGM}$ 
as quantum homogeneous spaces for a quantum $SL_{\infty}$. 
In the language of quantum groups
this means that we need to construct deformations of the two rings
$\widehat{k[\d_{\infty}]}$ and $k[\d_{\infty}]$ together with a
coaction of a deformation of ${k[SL_{\infty}]}$ on them.
The naturality of the construction in \S 2 will allow us to repeat the 
same arguments used for the commutative case also in the non commutative
case with very small changes.

{\bf Definition (3.1)}.
Let $k_q=k[q,q^{-1}]$ and
let $k_q<a_{i,j}>_{m,n}$ be the free algebra over $k_q$
with $a_{ij}$ as non commutative
generators, $-m \leq i, j \leq n-1$.
Define $k_q[a_{i,j}]_{m,n}$,
as the associative $k_q$-algebra with unit 
generated by the elements $a_{ij}$, subject
to the relations: \vskip1ex
$a_{ij}a_{kj}=q^{-1}a_{kj}a_{ij}, \quad i<k, \quad \quad 
a_{ij}a_{kl}=a_{kl}a_{ij},\quad i<k,j>l \quad or \quad i>k,j<l$ \vskip1ex
$a_{ij}a_{il}=q^{-1}a_{il}a_{ij}, \quad j<l,  \quad \quad 
a_{ij}a_{kl}-a_{kl}a_{ij}=(q^{-1}-q)a_{kj}a_{il}, \quad i<k,j<l$\vskip1ex

$k_q[a_{i,j}]_{m,n}$ is a bialgebra with counit and comultiplication:
$$
\epsilon^q_{(m,n)}(a_{ij})=\d_{ij} 
\quad \D^q_{(m,n)}(a_{ij})=\sum a_{ik} \otimes a_{kj}
$$
See [Ma1], [Ma2] for more details.

{\bf Definition (3.2)}.
We define the {\it quantum determinant} {\it obtained by taking 
rows \break $i_1 \dots i_p$,} {\it columns $j_1 \dots j_p$}
as an element $D_{i_1 \dots i_p}^{j_1 \dots j_p}$ $\in$ 
$k_q<a_{i,j}>_{m,n}$ given by:
$$
D_{i_1 \dots i_p}^{j_1 \dots j_p}=_{def}\sum_{\sigma:(i_1 \dots i_p)
\rightarrow (j_1 \dots j_p)}
(-q)^{-l(\sigma)}a_{i_1\sigma(i_1)} \dots a_{i_m\sigma(i_p)},
\quad \matrix{-m \leq i_1 < \dots < i_p \leq n-1 \cr
              -m \leq j_1 < \dots < j_p \leq n-1 }
$$
where $\sigma$ runs over all the bijections and $l(\sigma)$ is the
length of the permutation $\sigma$.
$p$ is called the {\it rank} of $D_{i_1 \dots i_p}^{j_1 \dots j_p}$.
Its image in $k_q[a_{i,j}]_{m,n}$ is then the usual quantum determinant.
We shall write $D_{i_1 \dots i_p}^{j_1 \dots j_p}$ for this image
also, the context making clear where the element sits. (See [PW] ch. 4
for more details). We will drop the upper indices whenever they
coincide with $-p \dots -1$.

{\bf Definition (3.3)}. Define {\it the quantum grassmannian ring} 
$k_q[\D_{m,n}]$,
as the subring of $k_q[a_{i,j}]_{m,n}$ generated by the 
quantum determinants $D_{i_0 \dots i_{m-1}}$
$-m \leq i_0 < \dots < i_{m-1} \leq n-1$ (see [Fi1]).
We will refer to the set of such determinants with $\D_{m,n}$.

An explicit presentation of the
ring $k_q[\D_{m,n}]$ in terms of generators and relations is given by 
(see [TT] 3.5, [Fi1], [FH]):
$$
\matrix{
q^{-[m-p]}\l_{J}\l_I = 
\l_I\l_{J}
+\sum_{i=1}^{N} (q^{-1}-q)^i
\cr \cr  
\sum_{(L,L') \in C^{i_0 \dots \hat i_{k_1} \dots \hat i_{k_p} \dots
i_{m-1}}_{i}} (-q)^{-l(\s(L))-l(\s(L'))}\l_{(L,i_{k_1} \dots i_{k_p})_{ord}}
\l_{(L',i_{k_1} \dots i_{k_p})_{ord}}
\cr \cr 
\quad I=(i_0 \dots i_{m-1}) < J=(j_0 \dots j_{m-1})
\quad I \cap J= \{i_{k_1} \dots i_{k_p}\}
\cr \cr
i_0<...<i_{m-1}, \quad j_0<...<j_{m-1}
\cr \cr
} \qquad (c)
$$
$$
\matrix{
\sum_{1 \leq \a_1 < ... < \a_s \leq m+s}
(-q)^{-l(z_1 \dots \hat z_{\a_1} ... \hat z_{\a_s} \dots
z_{r+s}z_{\a_1} ... z_{\a_s})-l(z_{\a_1} ... z_{\a_s}l_1 \dots l_{m-s})} 
\cr \cr
\l_{z_1 \dots \hat z_{\a_1} ... \hat z_{\a_s} \dots z_{m+s}}
\l_{(z_{\a_1} ... z_{\a_s}l_1 \dots l_{m-s})_{ord}}=0
\qquad (y)
}
$$
Each of the relations in the set $(y)$ is computed 
for any set of fixed indices: 
$-m\leq z_1 < \dots < z_{m+s}\leq n-1$, 
$-m \leq l_1 < \dots < l_{m-s} \leq n-1$.
\medskip
All the symbols that appear have been defined in [Fi2].
 
Notice that the relations labeled $(c)$ reduce for $q=1$ to state
the commutativity of the $\l_I$'s while the relations labeled 
$(y)$ for $q=1$ become the Young (also called symmetry) relations. 

We want now to proceed
in analogy with \S 2 and define the following inverse
and direct families.

{\bf Definition (3.4)}. 
Let $m'\geq m$, $n'\geq n$.
Define an inverse family or rings:
$$
\matrix
{ 
k_q[\D_{m',n'}] & \mapright {e^q_{(m',n',m,n)}} & k_q[\D_{m,n}] \cr\cr
D_{-m'\dots -m-1l_0 \dots l_{m-1}} & \longmapsto &
D_{l_0 \dots l_{m-1}}\cr\cr
D_{l_0 \dots l_{m'-1}} & \longmapsto &  
0 \ \ \hbox{otherwise}
} 
$$
where $-m \leq l_0 < ... < l_{m-1} \leq n-1$.

We define:
$$
\widehat {k_q[\D_{\infty}]}=\lim_{\longleftarrow} k_q[\D_{m,n}]
$$
Denote the induced maps
$$
e^q_{(m,n)}: \widehat {k_q[\D_{\infty}]}\lra k_q[\D_{m,n}].
$$
We observe that the maps $e_{(m',n',m,n)}$ are induced by maps
$$E^q_{(m',n',m,n)}:k_q[a_{i,j}]_{m',n'}\lra k_q[a_{i,j}]_{m,n}$$  
defined by
$$E^q_{(m',n',m,n)}(a_{ij})=a_{ij},\ \forall \ -m \leq i\leq n-1,\
-m\leq j\leq -1$$
$$E^q_{(m',n',m,n)}(a_{ij})=1,\ \forall \ -m'\leq i=j\leq -m-1$$
$$E^q_{(m',n',m,n)}(a_{ij})=0\ \hbox{otherwise}.$$
  
We define:
$$
k_q[M_{\infty }]=\lim_{\longleftarrow} k_q[a_{i,j}]_{m,n}
$$

{\bf Definition (3.5)}.
Let $m'\geq m$, $n'\geq n$.
Define the direct family of rings:

$$
\matrix
{
k_q[\D_{m,n}] & \mapright{r^q_{(m,n,m',n')}} & k_q[\D_{m',n'}] \cr\cr
D_{l_0 \dots l_{m-1}} & \longmapsto &
D_{-m', \dots -m-1,l_0 \dots l_{m-1}}
}
$$
for $-m \leq l_0 < \dots < l_{m-1} \leq n-1$.
We define: 
$$
k_q[\D_{\infty}]=\lim_{\longrightarrow} k_q[\D_{m,n}]
$$
Denote the induced inclusions $r^q_{(m,n)}:k_q[\D_{m,n}]\lra k_q[\D_{\infty}]$.

{\bf Observation (3.6)}.
Both $r^q_{(m,n,m',n')}$ and $e^q_{(m',n',m,n)}$ are well defined
that is they are zero on the relations on the determinants. This can
be directly checked.
\medskip
In analogy to \S 2, 
we can define a map:
$$
\matrix{
\matrix{
\rho^q_{(m,n)}:k[\D_{\infty}] & \longrightarrow & k[\D_{m,n}]
} \cr \cr
\rho^q_{(m,n)} (D_{a_{\bullet}})= 
\cases{
D_{a_{\leq m}} &  for all $m\geq ||a_{\bullet}||$, $n\geq |a_1|$ \cr \cr
0 &  otherwise
}
\cr
}
$$ 
Then we define
$$
\widehat {D}_{a_{\bullet}}=\{ \rho^q _{(m,n)}
D_{a_{\bullet}} \} \in \widehat{k[\D_{\infty}]}.
$$
\medskip

let $k_q<\xi _{a_{\bullet }}>$ to be the 
non commutative ring generated by the independent
variables $\xi _{a_{\bullet }}$, where ${a_{\bullet }}$ is any
Maya diagram of virtual cardinality 0. 

In analogy with \S 2
there is a natural map $\phi_q :k<\xi _{a_{\bullet }}>  
\lra {k[\D _{\infty }]}$
such that $\xi _{a_{\bullet }}  \lra  D_{a_{\bullet }}$.
This induces a topology on $k<\xi _{a_{\bullet }}>$ 
for which a fundamental set
of neighborhoods is given by $V^q_k:=\phi ^{-1}_qI^{-1}_qU_k^q$, where
$U^q_k$ is a fundamental set of neighbourhoods in 
$\widehat{k[\D _{\infty }]}$, defined by: 
$U^q_0=\widehat{k[\D _{\infty }]}$, 
$U^q_k=({e^q}_{(k,k)})^{-1}(0)$ and $I_q$ is the natural map from
$k[\D _{\infty }]$ and $k[\D _{\infty }]$ defined as 
$I_q(D _{a_{\bullet} })=\widehat D_{a_{bullet}}$. 

Let
$\widehat {k<\xi _{a_{\bullet }}>}$ be the completion of $
k<\xi _{a_{\bullet }}>$ with respect to the above
topology. In particular the elements of $\widehat {k<\xi _{a_{\bullet }}>}$
are of the form
$\sum b_i \xi _{a_{\bullet }^{i,1}}\cdots
\xi _{a_{\bullet }^{i,k_i}}$ where $b_i\in k$ and  the
${a_{\bullet}^{i,k}}$ are
any Maya diagrams of virtual cardinality 0.
The corresponding natural map between completions 
$\widehat {\phi^q } :\widehat {k<\xi _{a_{\bullet }}>} \lra 
\widehat {k[\D _{\infty }]}$ is defined by 
$\xi _{a_{\bullet }}   \lra \widehat{D}_{a_{\bullet }}$.


Define
$P_{(m,n),q}$ to be the two-sided ideal generated by
the relations (c) and (y) in $k_q[\D_{m,n}]$.
Let $P_q=\cup P_{(m,n),q}$ be the corresponding 
two-sided ideal in $k[\D_{\infty}]$,
and similarly $\widehat {P_q}=\matrix{\lim \cr {\longleftarrow}} 
P_{(m,n)}$
be the corresponding ideal in $\widehat {k_q[\D_{\infty}]}$.
As in the commutative case we have the following theorem.
\medskip

{\bf Theorem (3.7). Quantum deformation of the ind-variety 
$G_{\infty}$ and the Sato grassmannian $\widetilde{UGM}$.}

{\it 
$$
\matrix{
\widehat {k_q[\D _{\infty }]}\cong \widehat {k_q<\xi _{a_{\bullet }}>}/
\widehat {P_q}
\cr \cr
k_q[\D _{\infty }]\cong k_q<\xi _{a_{\bullet }}>/P_q.
}
$$
} 

{\bf Proof}. Same as (2.9).

{\bf Remark (3.8)}. If we specialize $q$ to 1 we obtain theorem (2.9),
that is a presentation of the commutative rings $k[\d_{\infty}]$
and $\widehat{k[\d_{\infty}]}$. Since by (2.11) 
$\hbox{Spec}_m({k[\d_{\infty}]})$ coincides with $G_{\infty}$
and $\hbox{Spec}_m(\widehat{k[\d_{\infty}]})$ with $\widetilde{UGM}$
we refer to the rings $k_q[\D_{\infty}]$
and $\widehat{k_q[\D_{\infty}]}$ as {\it quantum ind-variety $G_{\infty}$} 
and {\it quantum Sato grassmannian} respectively.
\medskip
Let's now proceed to show that the quantum ind-variety $G_{\infty}$
and the quantum Sato grassmannian are quantum homogeneous spaces.

{\bf Definition (3.9)}. In complete analogy with \S 2 define 
$$
k_q[SL_{m,n}]=k_q[a_{i,j}]_{m,n}/(D_{-m...n-1}^{-m...n-1}-1).
$$
$k_q[SL_{m,n}]$ is a quantum group, that is an Hopf algebra, with
antipode:
$$
S^q_{(m,n)}(a_{ij})=(-q)^{i-j}D^{-m...\hat j ... n-1}_{-m...\hat i ... n-1}.
$$
The coalgebra structure (i.e. the comultiplication and counit maps)
is naturally inherited from the matrix bialgebra $k_q[a_{i,j}]_{m,n}$.
For more details see [Ma1], [Ma2].

Then define the inverse system:
$$
\matrix{
k_q[SL_{m',n'}] \quad \mapright{\phi^q_{(m',n',m,n)}}  \quad k_q[SL_{m,n}]
\cr \cr
\phi_{(m',n',m,n)}^q(g_{ij})=
\cases{
g_{ij} & if $-m \leq i,j \leq n-1$
\cr \cr
1  & if $-m' \leq i=j \leq m'-1$, $n \leq i=j \leq n'-1$
\cr \cr
0  & otherwise
}
}
$$
One can directly check that these maps are well defined.

Define {\it quantum infinite dimensional special linear group}:
$$
{k_q[SL_{\infty}]}=_{def}\lim_{\longleftarrow} k_q[SL_{m,n}]
$$
Notice that for $q=1$ this coincides with the coordinate ring
of the ind-variety $SL_{\infty}(k)$.

We now intend to show that ${k_q[SL_{\infty}]}$ is a quantum group
that is it admits an Hopf algebra structure. This will be proved
in the same exact way we did for $k[SL_{\infty}]$ in \S 2. 

{\bf Proposition (3.10)}. {\it 
Let $(m',n')>(m,n)$. The following diagrams are commutative:

i)
$$
\matrix{
k_q[SL_{m',n'}] & \mapright{\phi^q_{(m',n',m,n)}} & k_q[SL_{m,n}] \cr \cr
\mapdown{\Delta^q_{(m',n')}} & & \mapdown{\Delta^q_{(m,n)}} \cr \cr
k_q[SL_{m',n'}] \otimes k_q[SL_{m',n'}] &
\mapright{\phi_{(m',n',m,n)} \otimes \phi^q_{(m',n',m,n)}} &
k_q[SL_{m,n}] \otimes k_q[SL_{m,n}]
}
$$

where $\Delta^q_{(m,n)}$ is the comultiplication
in $k_q[SL_{m,n}]$.

ii)
$$
\matrix{
k_q[SL_{m',n'}] & \mapright{\phi_{(m',n',m,n)}} & k_q[SL_{m,n}] \cr\cr 
\mapdown{\epsilon^q_{(m',n')}} & & \mapdown{\epsilon^q_{(m,n)}} \cr\cr
k_q & \mapright{id} & k_q
}
$$

where $\epsilon^q_{(m,n)}$ is the counit in $k_q[SL_{m,n}]$.

iii) 
$$
\matrix{
k_q[SL_{m',n'}] & \mapright{\phi_{(m',n',m,n)}} & k_q[SL_{m,n}] \cr \cr
\mapdown{S^q_{(m',n')}} & & \mapdown{S^q_{(m,n)}} \cr \cr
k_q[SL_{m',n'}] & \mapright{\phi_{(m',n',m,n)}} & k_q[SL_{m,n}]
}
$$

where $S_{(m,n)}$ is the antipode in $k_q[SL_{m,n}]$.
}

{\bf Proof}. Direct check. Notice that here the check involves also
the non commutative relations among the generators.

{\bf Corollary (3.11)}. {\it $k_q[SL_{\infty}]$ has an Hopf algebra
structure given by:

a) comultiplication
$$
\matrix{
k_q[SL_{\infty}] & \mapright{\Delta^q_{\infty}} & 
k_q[SL_{\infty}] \hat \otimes k_q[SL_{\infty}] \cr \cr
\{a_{(m,n)}\} & \mapsto & \{\Delta^q_{(m,n)}(a_{(m,n)})\}
}
$$

b) counit
$$
\matrix{
k_q[SL_{\infty}] &  \mapright{\epsilon^q_{\infty}} & k_q \cr \cr
\{a_{(m,n)}\} & \mapsto & \epsilon^q_{(m,n)}(a_{(m,n)})
}
$$

c) antipode
$$
\matrix{
k_q[SL_{\infty}] & \mapright{S^q_{\infty}} & k_q[SL_{\infty}]  \cr \cr
\{a_{(m,n)}\} & \mapsto & \{S^q_{(m,n)}(a_{(m,n)})\}
}
$$

where $\hat \otimes$ denotes the completed tensor product (see
[Ku] for more details).}
\medskip

Now we are ready to show that $\widehat{k[\D_{\infty}]}$ and 
$k[\D_{\infty}]$ are quantum homogeneous spaces.

We have the following coaction (see [Fi1]):
$$
\matrix{
k_q[\D_{m,n}] & \mapright{\lambda^q_{(m,n)}} & k_q[SL_{m,n}] \otimes k_q[\D_{m,n}] 
\cr \cr
D_{l_0 \dots l_{m-1}} & \longmapsto & 
\sum_{m \leq k_0 ... k_{m-1} \leq n-1} g_{l_0k_0}
...g_{l_{m-1}k_{m-1}} \otimes D_{k_0 ... k_{m-1}}
}
$$
One can check the commutativity of the following diagram, for
$m'\geq m$, $n'\geq n$:
$$
\matrix{
k_q[\D_{m,n}] & \mapright{\lambda^q_{(m,n)}} & k_q[SL_{m,n}] \otimes 
k_q[\D_{m,n}]
\cr \cr
\mapup{e^q_{(m',n',m,n)}} & & \mapup{\phi^q_{(m',n',m,n)} \otimes 
e^q_{(m',n',m,n)}} \cr \cr
k_q[\D_{m',n'}] & \mapright{\lambda^q_{(m',n')}} & k_q[SL_{m',n'}] 
\otimes k_q[\D_{m',n'}]
}
$$

{\bf Proposition (3.12)}. {\it 
There is an coaction of ${k_q[SL_{\infty}]}$ on
$\widehat{k_q[\D_{\infty}]}$ and on $k_q[\D_{\infty}]$.}

{\bf Proof}. Same as (2.18).
%
%
%
%
%
\bigskip
{\bf 4. Infinite dimensional invariant theory for $SL_{\infty}$ 
and its quantum deformation $k_q[SL_{\infty}]$}
\medskip
The first fundamental theorem of invariant theory for the special linear
group 
\break
$SL_{m,0}(k)$ (see (2.13) for the notation) 
states that given the right action of $SL_{m,0}(k)$ 
on the matrix algebra $k[b_{i,j}]_{m,n}$, where $m < n$,
$-m \leq i \leq n-1$, $-m \leq j \leq -1$:
$$
\matrix{
k[b_{i,j}]_{m,n} \times SL_{m,0}(k) & \lra & k[b_{i,j}]_{m,n} \cr \cr
(b_{ij},g) & \longmapsto & \sum b_{ik}g_{kj}
}
$$
the subring of invariants $k[b_{i,j}]_{m,n}^{SL_{m,0}(k)}$ coincides with the
subring generated by the determinants of rank $m$ in 
$k[b_{i,j}]_{m,n}$.
We want to generalize this result to the infinite dimensional case.

{\bf Observation (4.1)}. There is a natural right action of
$SL_{m,0}(k)$ on $M_{m,n}$ the set of matrices with row indices
from $-m$ to $n-1$ and column indices from $-1$ to $-m$.
$$
\matrix{
M_{m,n} \times SL_{m,0}(k) & \lra & M_{m,n}  \cr \cr
A,g & \longmapsto & Ag
}
$$
This action gives rise to the following coaction:
$$
\matrix{
k[b_{i,j}]_{m,n} & \mapright{\rho_{(m,n)}} 
& k[b_{i,j}]_{m,n} \otimes k[SL_{m,0}] \cr \cr
b_{ij} & \longmapsto & \sum b_{ik} \otimes g_{kj}
}
$$
An element $x \in k[b_{i,j}]_{m,n}$ is said to be coinvariant under this
coaction if $\rho_{(m,n)}(x)=x \otimes 1$. 
The first theorem of coinvariant theory equivalently states that
the subring of coinvariants under the coaction $\rho_{(m,n)}$,
$k[b_{i,j}]_{m,n}^{k[SL_{m,0}]}$ coincides with
$k[\d_{m,n}]$.

{\bf Proposition (4.2)}. {\it There is a coaction $\rho_{\infty}$
of $k[SL_{\infty ,0}]:=lim _{\leftarrow} k[SL_{m,0}]$ on $k[M_{\infty}]$.}

{\bf Proof.} Fix an index $m_0$ and for $m,n \geq m_0$ define the map:
$$
\matrix{
k[b_{i,j}]_{m,n} & \mapright{\rho_{(m,n,m_0)}}
& k[b_{i,j}]_{m,n} \otimes k[SL_{m_0,0}] \cr \cr
b_{ij} & \mapsto & \sum_{k=-m_0}^{n-1} b_{ik} \otimes g_{kj}+
\sum_{k=-m'}^{-m_0-1} b_{ik} \otimes \delta_{kj}g_{kj}
}
$$
where $\delta_{kj}=1$ if $k=j$ and 0 otherwise.

For any indices $m'\geq m\geq m_0$, $n'\geq n\geq 
n_0$ we have the commutative diagram:
$$
\matrix{
k[b_{i,j}]_{m,n} & \mapright{\rho_{(m,n,m_0)}} 
& k[b_{i,j}]_{m,n} \otimes k[SL_{m_0,0}] \cr \cr
\mapup{e_{(m',n',m,n)}} & & \mapup{e_{(m',n',m,n)} \otimes id} \cr \cr
k[b_{i,j}]_{m',n'} & \mapright{\rho_{(m',n',m_0)}}
& k[b_{i,j}]_{m',n'} \otimes k[SL_{m_0,0}] 
}
$$
This gives us a map:
$$
\matrix{
k[M_{\infty}] & \mapright{\rho_{m_0}} &
k[M_{\infty}] \otimes k[SL_{m_0,0}]
}
$$
Going to the inverse limit we obtain a map:
$$
\matrix{
k[M_{\infty}]  & \mapright{\rho_{\infty}} &
 k[M_{\infty}] \hat \otimes k[SL_{\infty ,0}]
}
$$
which is the required coaction. QED.
\medskip

We remark that the natural inclusion $k[SL_{\infty ,0}]\lra k[SL_{\infty}]$
is not an isomorphism, however there exist non canonical isomorphisms
between these two rings.

Let $k[M_{\infty}]^{k[SL_{\infty ,0}]}$ denote the subring of
$k[SL_{\infty ,0}]$-coinvariants, that is of those elements
$X$ such that $\rho_{\infty}(X)=X \otimes 1$. It is easy to see that
$x=\{ x_{(m,n)} \}\in k[M_\infty ]$ is $k[SL_{\infty ,0}]$-coinvariant
iff each $ x_{(m,n)}\in k[a_{i,j}]_{m,n}$ is $k[SL_{m,0}]$-coinvariant,
i.e. if $\rho _{m,n} (x_{m,n})=a_{m,n}\otimes 1$.
\medskip
{\bf Theorem (4.3). The first fundamental theorem of coinvariant theory for
$SL_{\infty ,0}(k)$.} 
$$
k[M_{\infty}]^{k[SL_{\infty ,0}]}=\widehat{k[\d_{\infty}]}
$$
{\bf Proof}.
The fact that $\widehat{k[\d_{\infty}]} \subset
k[M_{\infty}]^{k[SL_{\infty ,0}]}$ can be shown by checking directly
that the generators of $\widehat{k[\d_{\infty}]}$ are coinvariant,
that is:
$$
\rho_{\infty}(\widehat{d}_{a_{\bullet}})=
\widehat{d}_{a_{\bullet}} \otimes 1.
$$
For the other inclusion, let $x \in k[M_{\infty}]^{k[SL_{\infty ,0}]}$.
We need to prove that $x$ can be written as:
$$
x=\sum x_i \widehat{d}_{a_{\bullet}^{i,1}}
\dots \widehat{d}_{a_{\bullet}^{i,i_k}}
$$
This can be done using exactly the same argument as in (2.7)(b).
\medskip
We now turn to examine the quantum case.

In [FH] we prove that there is a well defined coaction:
$$
\matrix{
k_q[b_{i,j}]_{m,n} & \mapright{\rho_{m,n}^q} & 
k_q[b_{i,j}]_{m,n} \otimes k_q[SL_{m,0}]
\cr \cr
b_{ij} & \longmapsto & \sum b_{ik} \otimes g_{kj}
}
$$
and that:
$$
k_q[b_{i,j}]_{m,n}^{k_q[SL_{m,0}]}=k_q[\D_{m,n}]
$$ 

{\bf Theorem (4.4). The first fundamental theorem of quantum coinvariant
theory for $k_q[SL_{\infty ,0}]$}.

{\it There is a natural right coaction $\rho_{\infty}$ 
of $k_q[SL_{\infty ,0}]$ on
$k_q[M_{\infty}]$. Under this coaction the ring of coinvariants
coincides with the quantum infinite dimensional grassmannian
$\widehat{k_q[\D_{\infty}]}$ i.e.
$$
k_q[M_{\infty}]^{k_q[SL_{\infty ,0}]}=\widehat{k_q[\D_{\infty}]}.
$$
}

{\bf Proof.} Same as (4.2) and (4.3). 
\bigskip

\medskip
\centerline{\bf REFERENCES}
\medskip
\item{[A]} Abe, E., {\it Hopf Algebras}, Cambridge University Press,
Cambridge, MA, 1980.
\vskip1ex
\item{[CP]} Chari, V.; Pressley, A., {\it 
A guide to quantum groups}, 
Cambridge Press, 1994. 
\vskip1ex
\item{[FRT]} Faddeev, L. D.; Reshetikhin, N. Y.; Takhtajan, L. A., 
{\it Quantization of 
Lie groups and Lie algebras}, Algebraic analysis, Vol. I, 
129-139, Academic Press, Boston, MA, 1988.
\vskip1ex
\item{[Fi1]} Fioresi, R., {\it Quantization of the grassmannian manifold},
J. Algebra, {\bf 214}, no. 2, 418-447, 1999.
\item{[Fi2]} Fioresi, R., {\it Quantum deformation of the flag manifold}
Comm. Alg., {\bf 27}, no. 11, 5669-5685, 1999.
\vskip1ex
\item{[FH]} Fioresi, R.; Hacon C. {\it Quantum coinvariant theory
and quantum Schubert manifolds}, preprint, 1999.
\vskip1ex
\item{[Fu]} Fulton, W., {\it Young Tableaux}, Cambridge University Press,
Cambridge, 1997.
\vskip1ex
\item{[GH]} Griffiths, P.; Harris, J., 
{\it Principles of algebraic geometry},
 Wiley and Sons, New York, 1978. 
\vskip1ex
\item{[Ka]} Kac, V.; Raina, A. K., 
{\it Bombay Lecture notes on weight representations
of infinite dimensional Lie algebras}, Advanced Series in Mathematical
Physics, World Scientific Publishing, 1987.
\vskip1ex
\item{[Ku]} Kumar, S., {\it Infinite grassmannians and moduli spaces
of $G$-bundles}, Lecture Notes in Mathematics, 1649, 1-49, Springer, 1997.
\vskip1ex
\item{[La]} Lang, S. {\it Hilbert's Nullstellensatz in infinite
dimensional space}, Proc. Amer. Math. Soc., {\bf 3}, 407-410, 1952.
\vskip1ex 
\item{[LR]} Lakshmibai, V.; Reshetikhin, N., 
{\it Quantum flag and Schubert schemes},
(Amherst, MA, 1990), 145-181, Contemp. Math., {\bf 134}, Amer. Math. Soc., 
Providence, RI, 1992.
\vskip1ex
\item{[Ma1]} Manin, Y. I., {\it Quantum groups and non commutative geometry},
Centre de Reserches Mathematiques Montreal, 49, 1988. 
\vskip1ex
\item{[Ma2]} Manin, Y. I., {\it Topics in noncommutative geometry}, 
M. B. Porter Lectures,
Princeton University Press, Princeton, NJ, 1991.  
\vskip1ex
\item{[Mu1]} Mulase, M., {\it A correspondence between an infinite dimensional
grassmannian and arbitrary vector bundles on algebraic curves}, 
Proceedings of Symposia in Pure Mathematics, 49,
Part I, Amer. Math. Soc., Providence RI, 1989.
\vskip1ex
\item{[Mu1]} Mulase, M., {\it Category of vector bundles on algebraic curves
on infinite dimensional grassmannians}, Int. J. Math, {\bf 1}, no. 3,
293-342, 1990.
\vskip1ex
\item{[PW]} Parshall, B.; Wang, J. P., {\it Quantum linear groups}, 
Mem. Amer. Math. Soc. 89, no. 439, 1991. 
\vskip1ex
\item{[PM]} Plaza Martin, F., {\it Grassmannian of $k((z))$:
Picard group equations and automorphisms} preprint, 1998.
\vskip1ex
\item{[PS]} Pressley, A.; Segal, G., {\it Loop Groups}, 
Cambridge Press, 1991.
\vskip1ex
\item{[Sa1]} Sato, M., {\it Soliton equations as dynamical systems on an
infinite dimensional grassmannian manifold}, Kokyoroku, RIMS, Kyoto
University, {\bf 439}, 30-46, 1981.
\vskip1ex
\item{[Sa2]} Sato, M., {\it The KP hierarchy and infinite dimensional
grassmannian manifold}, Proceedings of Symposia in Pure Mathematics, 49,
Part I, Amer. Math. Soc., Providence RI, 1989.
\vskip1ex
\item{[SW]} Segal, G.; Wilson G., {\it Loop groups and equations of KdV
type}, Publ. Math. IHES, {\bf 61}, 5-65, 1985.
\vskip1ex
\item{[So]}  Soibelman, Y. S., {\it On the quantum flag manifold}, 
(Russian)
Funktsional. Anal. i Prilozhen. {\bf 26} 1992, no. 3, 90-92; 
translation in Functional Anal. Appl. {\bf 26}, no. 3, 225-227, 1992.
\vskip1ex
\item{[TT]} Taft, E.; Towber, J., {\it 
Quantum deformation of flag schemes and Grassmann
schemes, I. 
A $q$-deformation of the shape-algebra for ${\rm GL}(n)$},
J. Algebra {\bf 142},
no. 1, 1-36, 1991. 
\vskip1ex
\item{[VS]}  Vaksman, L. L.; Soibelman, Y. S., {\it 
On some problems in the theory of
quantum groups}, Representation theory and dynamical systems, 3-55, 
Adv. Soviet Math., 9,
Amer. Math. Soc., Providence, RI, 1992.
\vskip1ex
\item{[W]} Witten, E., {\it Quantum Field Theory, Grassmannians
and algebraic curves}, Comm. Math. Phys., {\bf 113}, 529-600, 1988.
\vskip1ex

\end